\documentclass{article} \usepackage{amssymb}

\def\Alg{{\mathfrak A}}

\begin{document}

\bibliographystyle{unsrt}

\parindent0pt

\begin{center} {\large \bf Factorizations in universal enveloping 
algebras of three dimensional Lie algebras and generalizations}

\end{center} \vspace{3mm}

\begin{center}  {\sc Stephen Berman}$^1$, {\sc Jun Morita}$^2$ and
{\sc Yoji Yoshii}$^3$

\vspace{5mm}

 1) Department of Mathematics, University of Saskatchewan\\
 Saskatoon,Saskatchewan, Canada S7N 5E6\\
berman@snoopy.usask.ca\\

 2) Institute of Mathematics, University of Tsukuba\\
 Tsukuba, Ibaraki, 305-8571, Japan\\
morita@math.tsukuba.ac.jp\\

 3) Department of Mathematics, University of Alberta\\
 Edmonton,Alberta, Canada, T6G 2G1\\
yoshii@math.ualberta.ca\\

\end{center}

\vspace{3mm}

\begin{center}{\bf This paper is dedicated to Robert V. Moody on the
occasion of
his 60th birthday}

\end{center}
\vspace{3mm}

\vspace{5mm}

\begin{abstract} We introduce the notion of Lie algebras with
plus-minus pairs as well as regular plus-minus pairs. These notions deal with
certain factorizations in universal enveloping algebras. We  show that many
important
Lie algebras have such pairs and we classify, and give a full treatment of, the
three dimensional Lie algebras with plus-minus pairs.

\vspace{5mm}
\par \noindent
{\bf AMS Classification Numbers 17B05, 17B35, 17B67, 17B70}
\end{abstract}

\parindent15pt \vspace{10mm}

\section{Introduction}

All of our algebras will be over a field $F$ of characteristic zero. We
begin by
recalling the well known fact that if $L$ is a  Kac-Moody Lie algebra with the
usual Chevalley 
generators $\{e_i, f_i |1 \leq i \leq l \}$ 
satisfying
$L = \left[ L,L \right]$ and $l$ is fnite,  then every 
$L$-module on which the elements $e_i,f_i, 1 \leq i \leq l$ act locally
nilpotently
 is integrable in the sense that the elements $h_i=[e_i,f_i], 1\leq i \leq l$
are simultaneously diagonalizable. (cf. \cite{MP} Ex. 6.31,
p. 585, or \cite{Tits}).
In other words, every weakly integrable module for such an algebra is
integrable. The usual proof of this fact uses that the three dimensional Lie
algebra with basis  $e_i,f_i, h_i$ is isomorphic to the Lie algebra 
$sl_2$ (so $\mathfrak g$ is generated by $sl_2$-triples) together with the
result which says that if $V$ is any module for $sl_2$ on which the standard
generators
$e,f$  of $sl_2$ act locally nilpotently then the element $h=[e,f]$ is
diagonalizable on $V$. 
One can see this last fact as follows. We denote by $M(W)$
the maximal integrable submodule of an $sl_2$-module $W$. In general,
$M(W/M(W)) = 0$ for all $sl_2$-modules $W$. If a vector $v$ of an 
$sl_2$-module $W$ 
satisfies that $fv = 0$, $e^nv \not= 0$ and $e^{n+1}v = 0$, then
we obtain $h(e^nv) = n(e^nv)$, since $fe^n = e^nf - ne^{n-1}(h + n - 1)$ and
$$\begin{array}{lll}
h(e^nv) & = & (ef-fe)e^nv\\
& = & efe^nv\\
& = & e(e^nf - ne^{n-1}(h+n-1))v\\
& = & -ne^n(h+n-1)v\\
& = & -n(h-n-1)e^nv\\
& = & -nh(e^nv) + n(n+1)e^nv.
\end{array}$$
This implies that $M(V)$ is nontrivial
for every nonzero $sl_2$-module $V$ on which the elements $e$ and $f$ are
locally nilpotent operators. 
Therefore, $M(V) = V$ for such an $sl_2$-module $V$,
that is, $V$ is integrable.
\bigskip

We want to indicate another approach to the above fact about $sl_2$-modules
which uses a factorization in the universal enveloping algebra. This method
appears to be new and was the starting point of this paper. For any Lie algebra
$\mathfrak g$ we let $U(\mathfrak g)$ denote its universal enveloping algebra.
Then one knows that for the algebra $sl_2 =Fe \oplus Fh \oplus Ff$
we have the factorization
\begin{eqnarray}\label{sl2f}
U(sl_2)=U(Fe)U(Ff)U(Fe).
\end{eqnarray}
Using this it is easy to see that if $V$ is an $sl_2$-module on which both $e$
and $f$ act locally nilpotently then any vector $v \in V$  generates a finite
dimensional submodule. Thus $h$ acts semisimply on this submodule and so
 we obtain $h$ acts semisimply on $V$. Also, the proof of the factorization 
(\ref{sl2f})
is quite straightforward and follows easily from  the following formula in
$U(sl_2)$.
\begin{eqnarray}\label{mf}
f(e^if^je^k)=\frac{j-i+1}{j+1}e^if^{j+1}e^k +
        \frac{i}{j+1}e^{i-1}f^{j+1}e^{k+1} + i(j-i+1)e^{i-1}f^je^k, 
\end{eqnarray}
for all $i>0, j ,k \geq 0.$
This formula can be established using the following: for any $k \geq 0$ we have
$$\begin{array}{llll}
({\rm A}_k)\quad & fef^k & = & \frac{k}{k+1}ef^{k+1} + \frac{1}{k+1}f^{k+1}e
+ kf^k\\
({\rm B}_k)\quad & f^kef & = & \frac{1}{k+1}ef^{k+1} + \frac{k}{k+1}f^{k+1}e
+ kf^k.
\end{array}$$
$fe^i=e^if-ie^{i-1}(h+i-1)$ for $i \geq 1$ and induction.

Next let ${\mathfrak H}$ be the three dimensional Heisenberg Lie algebra with a
basis
$x,y,z$ satisfying $[x,y] = z,\ [x,z] = [y,z] = 0$. Then,
in $U({\mathfrak H})$, we obtain the following factorization
\begin{eqnarray}\label{heif}
U(\mathfrak H)=U(Fx)U(Fy)U(Fx).
\end{eqnarray}
The proof of this is much like the $sl_2$ case. It follows easily from the
following formula in $U(\mathfrak H)$
\begin{eqnarray}\label{heimf}
y(x^iy^jx^k) = \frac{j-i+1}{j+1} x^iy^{j+1}x^k
+ \frac{i}{j+1}x^{i-1}y^{j+1}x^{k+1},
\end{eqnarray}
for all $i>0, j ,k \geq 0.$
Note that (\ref{heimf}) is proved by establishing for any $k \geq 0$ we have
$$\begin{array}{llll}
({\rm A}_k')\quad & yxy^k & = & \frac{k}{k+1}xy^{k+1} + \frac{1}{k+1}y^{k+1}x\\
({\rm B}_k')\quad & y^kxy & = & \frac{1}{k+1}xy^{k+1} + \frac{k}{k+1}y^{k+1}x,
\end{array}$$
which in turn is proved using $yx^i=x^iy-ix^{i-1}z$ for $i \geq 1$ and
induction.
Thus, the picture is similar for both algebras $sl_2$ and $\mathfrak H$ in that
they both have a pair of subalgebras
$P,M$ satisfying $P + M$ is not the whole algebra and $U(P)U(M)U(P)$ is the
whole enveloping algebra. This prompts the following definition which singles
out those Lie algebras having this type of factorization in their universal
enveloping algebras.
\bigskip

{\bf Definition 1.1} (i.) A Lie algebra $L$ is said to have a  {\bf plus-minus
pair} if it has two subalgebras $P,M$ satisfying $P+M \neq L$ and
$$
U(L)=U(P)U(M)U(P).
$$
In this case we say $L$ has a plus-minus pair $(P,M)$.

(ii.) Let $(P,M)$ be a plus-minus pair of $L$. We say this is a {\bf regular
plus-minus pair} if 
$P \cap M= (0)$ and there is an automorphism $\sigma$ of $L$ of order two
satisfying $\sigma (P)=M.$ Note that in this case we then have
$U(L)=U(P)U(M)U(P)=U(M)U(P)U(M).$
\bigskip

It is clear that both Lie algebras $sl_2$ and $\mathfrak H$ have regular
plus-minus pairs. Moreover if $L$ is any three dimensional Lie algebra with a
plus-minus pair $(P,M)$ then each of $P$ and $M$ must be one dimensional.
Indeed, $P+M$ 
cannot be $3$ dimensional as $P+M \neq L$. Thus, $P+M$ is two dimensional
and so
if one of $P,M$ is $2$ dimensional then $P+M$ is a subalgebra of $L$ and so
$U(P+M) \neq U(L)$ but $ U(P)U(M)U(P) \subseteq U(P+M)$ so $(P,M)$ cannot be a
plus-minus pair. Letting 
$P=Fx, M=Fy$ we have that
\begin{eqnarray}\label{tpm}
U(L)= \sum_{i,j,k \geq 0} Fx^i y^j x^k.
\end{eqnarray}
Moreover,  the following result, which extends the situation discussed in the
$sl_2$ case, is quite clear.
 Let $L$ be a three dimensional Lie algebra with a plus-minus pair $(P,M)$
where

$P=Fx, M=Fy.$ Let $V$ be any $L$-module on which the action of the elements
$x,y$ is locally finite. Then any 
finitely generated submodule of $V$ is finite dimensional. Thus, one is led to
ask just which three dimensional Lie algebras have plus-minus pairs. 
\bigskip

   In Section 2 we will extend the methods used in the proofs for the
$sl_2$ and
$\mathfrak H$ cases above and show that any three dimensional Lie algebra which
is generated by two elements has a plus-minus pair. Then, using results on
three dimensional Lie algebras from \cite{Jac} we will see that there are only
two isomorphism classes of the three dimensional Lie algebras which do not have
plus-minus pairs. We also go on to study, when the base field $F$ is
algebraically closed, which of these algebras have regular plus-minus pairs and
are able to give a complete list of these. In the third and final section of
this paper we go on to investigate plus-minus pairs, or similar factorizations,
in the universal enveloping algebras, of Borcherds Lie algebras as well as in
some ${\bf Z}^n$-graded  Lie algebras which satisfy some extra conditions.

\section{Three dimensional case}\par

In this section we begin by showing a three dimensional Lie algebra has a
plus-minus pair if and 
only if it is generated by two elements. We then go on to investigate some
special cases as well as regular 
plus-minus pairs when the base field is algebraically closed.
\bigskip

   Throughout we let $L$ be a three dimensional Lie algebra unless mentioned
otherwise. If $(P,M)$ is a plus-minus pair for $L$ then we know that each of
$P,M$ is one dimensional so we let $P=Fx,M=Fy.$ If $x$ and $y$ do not generate
$L$ then it must be that $P+M$ is a proper subalgebra of $L$ and so since
$U(P)U(M)U(P) \subseteq U(P+M)$ we get a contradiction. Thus $L$ is
generated by
$x$ and $y$ so is two-generated.
We want to establish the converse of the above result. For the moment we
let $L$
be an arbitrary Lie algebra and 
$x, y$ any two elements of $L$.
\bigskip
  
We define subspaces $U_k$ for $k \geq 0$ of $U(L)$ by saying 
\begin{eqnarray}\label{Uk}
U_k = \sum_{0 \leq m \leq k}(Fxy^m +Fy^mx).
\end{eqnarray}
\medskip

Notice that $U_0=Fx $ and that $U_k \subseteq U_{k+1}$ for all $k \geq 0.$  
\bigskip

{\bf Lemma 2.1} Let $L$ be an arbitrary Lie algebra and $x,y$ any two elements
of $L$. 
 For $k \geq 0$ the following statements hold,
$$\begin{array}{ll}
({\rm A}_k) & yxy^k \equiv \frac{k}{k+1} xy^{k+1} + \frac{1}{k+1} y^{k+1}x
\quad {\rm mod}\ U_k,\\
({\rm B}_k) & y^kxy \equiv \frac{1}{k+1} xy^{k+1} + \frac{k}{k+1} y^{k+1}x
\quad {\rm mod}\ U_k,\\
({\rm C}_k) & y\ U_k \subset U_{k+1},\quad U_k\ y \subset U_{k+1}.
\end{array}$$
\bigskip

{\it Proof }.\ We prove this by induction on $k$ noting that for $k=0$ both
 $({\rm A}_0)$ and $({\rm B}_0)$ are clear. Next, we show $({\rm A}_0), \ldots,
({\rm A}_k), ({\rm B}_0), \ldots, ({\rm B}_k)$
imply $({\rm C}_k)$. Indeed, by definition we have that 
$$y\ U_k = \sum_{0 \leq m \leq k} (Fyxy^m + Fy^{m+1}x)$$ and so by $({\rm
A}_0),
\ldots, ({\rm A}_k)$ we get that this 
is contained in $U_{k+1}.$ Similiarly we have 
$$U_k\ y  =  \sum_{0 \leq m \leq k} (Fxy^{m+1} + Fy^mxy)$$ and so by $ ({\rm
B}_0), \ldots, ({\rm B}_k)$ we get this is contained in $U_{k+1}$. Hence $({\rm
C}_k)$ holds.
\medskip

Next we show that $({\rm A}_k), ({\rm B}_k), ({\rm C}_k)$ imply $({\rm
A}_{k+1}),
({\rm B}_{k+1})$. Now $yxy^{k+1} =(yxy^k)y$ so that $({\rm A}_k)$ implies that
the difference 
$$yxy^{k+1} -(\frac{k}{k+1} xy^{k+1} + \frac{1}{k+1} y^{k+1}x)y \in U_k y.$$
But $({\rm C}_k)$ implies that $U_k y \subseteq U_{k+1}$ so we get that
$$yxy^{k+1} \equiv (\frac{k}{k+1} xy^{k+1} + \frac{1}{k+1} y^{k+1}x)y  
\quad {\rm mod}\ U_{k+1}.$$ 
Similarly, using $({\rm B}_k)$ and $({\rm C}_k)$ we get that
$$y^{k+1}xy \equiv y
(\frac{1}{k+1} xy^{k+1} + \frac{k}{k+1} y^{k+1}x)\quad {\rm mod}\ U_{k+1}.$$
Thus, we finally get
$$yxy^{k+1} \equiv \frac{k}{k+1} xy^{k+2} + \frac{1}{(k+1)^2} yxy^{k+1} +
\frac{k}{(k+1)^2} y^{k+2}x \quad {\rm mod}\ U_{k+1}$$
which implies that
$$\frac{k(k+2)}{(k+1)^2} yxy^{k+1} \equiv
\frac{k}{k+1} xy^{k+2} + \frac{k}{(k+1)^2} y^{k+2}x\quad {\rm mod}\ U_{k+1}.$$
Therefore, we obtain
$$yxy^{k+1} \equiv
\frac{k+1}{k+2} xy^{k+2} + \frac{1}{k+2} y^{k+2}x\quad {\rm mod}\ U_{k+1},$$
and we see that $({\rm A}_{k+1})$ holds.
\medskip

Using a similiar type of argument we have that
$$\begin{array}{llll}
y^{k+1}xy & = & y(y^kxy) &\\
& \equiv & y(\frac{1}{k+1} xy^{k+1} + \frac{k}{k+1} y^{k+1}x) &
{\rm mod}\ U_{k+1}\\
& \equiv & \frac{1}{k+1} yxy^{k+1} + \frac{k}{k+1} y^{k+2}x &
{\rm mod}\ U_{k+1}\\
& \equiv & \frac{1}{k+1} (\frac{k}{k+1} xy^{k+1} +
\frac{1}{k+1} y^{k+1}x)y + \frac{k}{k+1} y^{k+2}x & {\rm mod}\ U_{k+1}\\
& \equiv & \frac{k}{(k+1)^2} xy^{k+2} + \frac{1}{(k+1)^2} y^{k+1}xy +
\frac{k}{k+1} y^{k+2}x & {\rm mod}\ U_{k+1}
\end{array}$$
and
$$\frac{k(k+2)}{(k+1)^2} y^{k+1}xy \equiv
\frac{k}{(k+1)^2} xy^{k+2} + \frac{k}{k+1} y^{k+2}x\quad {\rm mod}\ U_{k+1}.$$
Therefore, we obtain
$$y^{k+1}xy \equiv
\frac{1}{k+2} xy^{k+2} + \frac{k+1}{k+2} y^{k+2}x\quad {\rm mod}\ U_{k+1},$$
and we see that 
$({\rm B}_{k+1})$ holds. This completes our induction. Q.E.D.
\bigskip

We apply this Lemma to the three dimensional case in our next result.
\bigskip

{\bf Theorem 2.2} Let $L$ be a three dimensional Lie algebra. Then $L$ has a
plus-minus pair if
 and only if $L$ is two generated. Moreover, if $x$ and $y$ generate $L$ then
$(P,M)$ is a plus-minus pair for $L$ where 
$P=Fx,M=Fy.$
\bigskip

{\it Proof}. We need only show $L$ has a plus-minus pair if $L$ is generated by
two elements 
$x,y$. Let $z=[x,y].$ Now we want to show
$U(L) = \sum_{i,j,k \geq 0} Fx^iy^jx^k$. Put
${\mathfrak X} =
\sum_{i,j,k \geq 0} Fx^iy^jx^k \subset U(L)$ and let $U_k$ be defined as above.
Clearly $x {\mathfrak X} 
\subset {\mathfrak X},  {\mathfrak X}x 
\subset {\mathfrak X}$ and $U_k \subseteq {\mathfrak X}$ for all $k \geq 0$. We
claim
$$\begin{array}{lll}
y(x^\ell y^m x^n) \in {\mathfrak X},\\
z(x^\ell y^m x^n) \in {\mathfrak X}.
\end{array}$$
and show this by induction on $\ell$. If $\ell = 0$, then we see
$y(y^mx^n) \in {\mathfrak X}$ and using $({\rm A}_m)$ we get
$$\begin{array}{lll}
z(y^mx^n) & = & (xy - yx)(y^mx^n)\\
& = & xy^{m+1}x^n - yxy^mx^n\\
& \in & Fxy^{m+1}x^n + (Fxy^{m+1} + Fy^{m+1}x + U_m)x^n \subset  {\mathfrak X}.
 
\end{array}$$
Let $\ell > 0$. Then, we obtain, using our inductive assumption, that
$$\begin{array}{lll}
y(x^\ell y^m x^n) & = & (xy -z)(x^{\ell -1}y^mx^n)\\
& \in & x{\mathfrak X} + {\mathfrak X} \subset {\mathfrak X}
\end{array}$$
and, letting $[z, x]=ax+by+cz$ for $a,b,c \in F$, we also get using our
inductive assumption that
$$\begin{array}{lll}
z(x^\ell y^mx^n) & = & (xz + ax + by + cz)(x^{\ell -1} y^m x^n)\\
& \in & x{\mathfrak X} + {\mathfrak X} + {\mathfrak X} + {\mathfrak X} 
\subset {\mathfrak X}.
\end{array}$$
Hence, $y {\mathfrak X} \subset {\mathfrak X}$. Since
${\mathfrak X}$ is a left ideal of $U(L)$ containing $1$, we obtain
${\mathfrak X} = U(L)$.
Therefore, $(P,M)$ with $P = Fx$ and $M = Fy$ is a plus-minus pair for
$L$. Q.E.D.
\bigskip

If our three dimensional Lie algebra $L$ is abelian it clearly does not have a
plus-minus pair. Also, we 
let $\mathfrak g$ be the three dimensional Lie algebra with basis $x,y, z$
satisfying
$$[x,y]=0,[x,z]=x, [y,z]=y.$$
Then for any elements $a,b,c, \alpha, \beta, \gamma \in F$ we have the very
special identity 
$$[ax+by+cz,\alpha x +\beta y +\gamma z]= \gamma(ax+by +cz)-c(\alpha x +\beta y
+\gamma z).$$
This clearly implies that $\mathfrak g$ is not two generated so does not have a
plus-minus pair.
Our next result shows that these are the only two kinds of
three dimensional Lie algebras which do not have plus-minus pairs.
In the proof we will make use of the classification of three dimensional
algebras presented in Jacobson's book,
\cite{Jac}, pg.11-14. We break our argument into 5 subcases (labled (a.)
through
(e.)) just as he does. Here
$L^{\prime}$ denotes the derived algebra of $L$.
\bigskip

{\bf Theorem 2.3} Let $L$ be a three dimensional Lie algebra. Then $L$ has a
plus-minus pair if and only if $L$ is  not abelian and is not isomorphic to the
algebra $\mathfrak g$ above.
\bigskip

{\it Proof}. (a.) If $L^{\prime} =(0)$ then $L$ is abelian and it clearly
has no
plus-minus pair.

(b.) Assume $L^{\prime}$ is one dimensional and in the center, $Z(L)$, of $L$.
Then there is only one Lie algebra staisfying these conditions and it has basis
$x,y,z$ satisfying $[x,y] = z, [x,z] = [y,z] = 0$. This is clearly generated by
$x,y$ so has a plus-minus pair.

(c.) Assume that $L^{\prime}$ is one dimensional but not in $Z(L)$. Then
$L$ has
basis $x,y,z$ satisfying $[x,y] = x, [x,z] = [y,z] = 0$ so $L$ is generated by
the two elements $x+z$ and $y$, hence has a plus-minus pair.

(d.) Suppose that $L^{\prime}$ is two dimensional. Then, from \cite{Jac}, one
has that $L^{\prime}$ is abelian with basis, say $x,y$, and there is an element
$z \in L$ so that  $[x,y] = 0, [x,z] = \alpha x + \beta y,
[y,z] = \gamma x + \delta y$ with $\alpha \delta - \beta \gamma \not= 0$.
If $\alpha = 0$, then $x$ and $y + z$ generate $L$ since $[x,y+z] = \beta
y$ and
$\beta \neq 0$ as $\alpha =0.$

Suppose $\alpha \not= 0$. By
replacing $z$ by $z/\alpha$, we can assume $\alpha = 1$. Hence,
we have that
$[x,y] = 0, [x,z] = x + \beta y, [y,z] = \gamma x + \delta y$ with
$\delta \not= \beta \gamma$. If $\gamma \not= 0$, then
$- \gamma x + y$ and $z$ generate $L$ since
$[- \gamma x + y , z] = (\delta - \beta \gamma)y$.
Next suppose $\gamma = 0$. Then 
we have $\delta \not= 0$, and our relations are $[x,y] = 0, 
[x,z] = x + \beta y, [y,z] = \delta y$. If $\beta \not= 0$, then
$x - (\beta/\delta)y$ and $z$ generate $L$ since
$[x - (\beta/\delta) y , z] = x$. In the remaining case, our relations are
$[x,y] = 0, [x,z] = x, [y,z] = \delta y$ with $\delta \not= 0$. If
$\delta \not= 1$, then $x+y$ and $z$ generate $L$ since  $[x+y,z] = x + \delta
y$. Therefore the only case left is when 
$\delta =1$ and then our algebra $L$ is nothing but the algebra $\mathfrak g$
above.

(e.) This is the case where $L^{\prime} =L$. Then Jacobson shows there is a
basis $x,y, z$
 of $L$ satisfying  $[x,y] = z, [x,z] = \alpha y, [y,z] = \beta x$ for some 
$\alpha, \beta \in F$ with $\alpha \beta \not= 0$. Here we see that $x$ and $y$
generate $L$. Q.E.D.
\bigskip

The special case when 
$L = L_{-1} \oplus L_0 
\oplus L_1$ is a three graded
Lie algebra of dimension three with a plus-minus pair can be discussed without
reference to Jacobson's classification.
We shall use this in the final section of this work so
we briefly go through this here. Put $L_1 = Fx,\ L_{-1} = Fy,\ 
L_0 = Fz$. We can assume first that
$[x,y]$ is either $0$ or $z$. Suppose $[x,y] = 0$. If $[x,z] = [y,z] = 0$,
then $L$ is abelian so has no plus-minus pair. If
$[x,z] = 0$ and $[y,z] \not= 0$, then we can also suppose $[z,y] = y$ and
hence, $P = F(x+y)$ and $M = Fz$ give a plus minus pair.
If $[x,z] \not= 0$ and $[y,z] = 0$, then we can suppose $[z,x] = x$ and
hence, again $P = F(x+y)$ and $M = Fz$ becomes a plus-minus pair.
If $[x,z] = ax$ and $[y,z] = by$ with $ab \not= 0$, then we can suppose
$a = 1$. In this case, $P = F(x+y)$ and $M = Fz$ give a plus-minus pair
when $b \not= 1$. Otherwise we have $[x,z] = x$ and $[y,z] = y$  and there
is no

plus-minus pair. Next we suppose $[x,y] = z$. If $[x,z] = [y,z] = 0$,
then $L$ is a Heisenberg Lie algebra, and hence, 
$L$ has a plus-minus pair. If $[x,z] = ax$ and $[y,z] = by$
with $a \not= 0$ or $b \not= 0$, 
then $0 = [z,z] = [[x,y],z] = [[x,z],y] + [x,[y,z]]
= a [x,y] + b [x,y] = (a+b)z$ and $a + b = 0$. Put $x' = x,\ y' = -2y/a,\ 
z' = -2z/a$. Then, $[x',y'] = -2[x,y]/a = -2z/a = z'$, $[z',x'] = -[x',z']
= 2[x,z]/a = 2x = 2x'$ and $[z',y'] = -[y',z'] = -4[y,z]/(a^2) = 
4y/a = -2y'$. This means that $L$ is isomorphic to $ sl_2$.
Therefore we obtain the following result which gives a
 characterization of $sl_2$ and $\mathfrak H.$ 
\bigskip

{\bf Proposition 2.4}.\ Let $L = L_1 \oplus
L_0 \oplus L_{-1}$ be a three graded Lie algebra
of dimension three with ${\rm dim}\ L_{\pm 1} = 
{\rm dim}\ L_0 = 1$.\par
\noindent
(1)\ If $L$ has a plus-minus pair, then $L$ is
isomorphic to one of $s\ell_2$, ${\mathfrak H}$ and $K(a,b)$, where
$K(a,b) = Fx \oplus Fy \oplus Fz$ is the Lie algebra having the 
relations:
$[x,y] = 0,\ [x,z] = ax,\ [y,z] = by$ with $a \not= b$.
\par
\noindent
(2)\ If $L$ has 
$(L_1,L_{-1})$ for a plus-minus pair, then $L$ is isomorphic to
either $sl_2$ or ${\mathfrak H}$.
\bigskip

{\bf Remark} If $a=b$ is non-zero then we have $K(a,b)=K(a,a)\simeq K(1,1)$ and
this is nothing but our algebra $\mathfrak g$ of Theorem 2.3 which does not
have
a plus-minus pair.
\bigskip

Next we will  briefly discuss isomorphism classes among the Lie algebras
$K(a,b)$.
Suppose that $K(a,b)$ is isomorphic to $K(a',b')$.
If $a$ or $b$ is $0$, then, by considering the derived subalgebra, we see that
one of $a'$ and $b'$ is $0$.
Furthermore we see $K(0,c) \simeq K(c,0) \simeq
K(0,1)$ for nonzero $c \in F$.
Next we suppose that both $a,b$ are nonzero.
Then, we also see $K(a,b) \simeq K(a/b,1)$.
Suppose $\phi : K(a,1) \rightarrow K(a',1)$
 is an isomorphism for nonzero $a,a' \in F$, where
$K(a,1)$ is generated by $x,y,z$ with $[x,y] = 0,\ 
[x,z] = ax,\ [y,z] = y$, and $K(a',1)$ is generated
by $x',y',z'$ with $[x',y'] = 0,\ [x',z'] = a'x',\ 
[y',z'] = y'$. Put $z'' = \phi(z) = a''x' + b''y' + c''z' \in
K(a',1)$. Then considering the derived algebras we find  $\phi(Fx \oplus Fy) =
Fx' \oplus Fy'$,
 and so we see that $c''$ is nonzero and that
$${\rm ad}\ z \mid_{Fx \oplus Fy} = 
\left( \begin{array}{cc} -a & 0\\ 0 & -1 \end{array} \right)$$
must be similar to
$${\rm ad}\ z'' \mid_{Fx' \oplus Fy'} =
\left( \begin{array}{cc} -c''a' & 0\\ 0 & -c'' \end{array} \right).$$
Therefore,
$a = c''a',\ 1 = c''$
or
$a = c'',\ 1 = c''a'$.
This means $a' = a^{\pm 1}$. Conversely, we can easily confirm
$K(a,1) \simeq K(1/a,1)$. We set ${\mathfrak P}(F) =
\{\ \{ u,u^{-1} \}\ \mid\ u \in F,\ u \not= 0,1\ \} \cup \{\ \{ 0 \}\ \}$. 
Then, the isomorphism classes of the Lie algebras $K(a,b)$,
having plus-minus pairs, are parametrized by
the set ${\mathfrak P}(F)$. Here the isomorphism class of $K(0,a)$ corresponds
to $\{ 0 \}$ while that of 
$K(a,1)$ to $\{ a,a^{-1} \}=\{ a^{-1}, a \}$ for $a \in F, a \neq 0.$ 
\bigskip

 We next assume that $F$ is an algebraically closed field
of characteristic $0$, and
 will study the three dimensional Lie algebras over $F$
having a regular plus-minus
pair. Let $L$ be such an algebra and let $(P,M)$ be a regular plus-minus
pair of
$L$. Then we can choose nonzero elements
$x \in P$ and $y \in M$ as well as an involutive automorphism $\sigma$
of $L$ such that $[x,y] \not= 0$ and $\sigma(x) = y$.
Put $z = [x,y]$, and set $u = x + y$ and $v = x - y$. Let
$L_1 = Fu$ (the $1$-eigenspace of $\sigma$)
and $L_{-1} = Fv \oplus Fz$ (the $-1$-eigenspace of $\sigma$).
Then, $[u,v] = [x+y,x-y] = -2 [x,y] = -2z$.
We write $[z,x] = ax + by + cz$. Then we obtain
$[z,y] = - \sigma([z,x]) = -\sigma (ax + by + cz) =
-bx -ay + cz$. The Jacobi identity implies
$$[[x,y],z] + [[y,z],x] + [[z,x],y] = 0$$ and this 
gives $[bx+ay-cz,x] + [ax+by+cz,y] = 0$. Therefore,
$$\{ -az - c(ax+by+cz) \} + \{ az + c(-bx-ay+cz) \} = 0$$
and $c(a+b)x + c(a+b)y = 0$, which implies
$c = 0$ or $a+b=0$.
\bigskip

\noindent
{\it Case 1}:\ $c = 0$.\par
\noindent
In this case, we have
$[x,y] = z,\ [z,x] = ax + by,\ [z,y] = - bx - ay$. We write the matrix of $\rm
ad (z)$ restricted to the space 
$Fx \oplus Fy$ as
$${\rm ad}\ z = \left(
\begin{array}{rr} a & -b\\ b & -a \end{array}
\right).$$
 Then, its characteristic polynomial 
is $t^2 - a^2 + b^2$. Hence,
${\rm ad}\ z\mid_{Fx \oplus Fy}$ is similar to one of
$$\left( \begin{array}{rr} \lambda & 0\\ 0 & - \lambda \end{array}
\right)\qquad \mbox{and}\qquad  
\left( \begin{array}{rr} 0 & 1\\ 0 & 0 \end{array}
\right),$$
where $\lambda = (a^2-b^2)^{1/2}$.
If $\lambda \not= 0$, that is $a^2 - b^2 \not= 0$, then
we have certain elements $x',y' \in Fx \oplus Fy$ such that
$[z',x'] = \lambda' x', [z',y'] = - \lambda' y'$ with
$z' = [x',y'] \not= 0$ and $\lambda' \not= 0$.
This means that $L \simeq sl_2$. If $a = b = 0$, then
we see $L \simeq {\mathfrak H}$. Now we suppose $a = b \not= 0$.
Thus, we have $[z,u] = 0, [z,v] = 2au, [u,v] = -2z$. Put
$\mu = (-a)^{1/4}$, and set $z' = z/\mu, u' = \mu u, v'= v/(2\mu^2)$.
Then,
$$\begin{array}{l}
\left[ v',u' \right] = [v,u]/(2\mu) = z/\mu = z',\\
\left[ v',z' \right] = (-a)u/(\mu^3) = u',\\
\left[ u',z' \right] = 0.
\end{array}$$
As is easy to check, one possible realization for this Lie algebra is obtained
by taking
$$v' = \left( \begin{array}{rrr} 
0 & 1 & 0\\ 1 & 0 & 0\\ 0 & 0 & 0
\end{array} \right),\qquad
u' = \left( \begin{array}{rrr}
0 & 0 & 1\\ 0 & 0 & 0\\ 0 & 0 & 0
\end{array} \right),\qquad
z' = \left( \begin{array}{rrr}
0 & 0 & 0\\ 0 & 0 & 1\\ 0 & 0 & 0
\end{array} \right),$$
and
$$\sigma(v')=-v',\sigma(u')=u', \sigma(z')=-z'.$$
\smallskip

Next we suppose $a = -b \not= 0$. Thus, we have
$[z,v] = 0, [z,u] = 2av, [u,v] = - 2z$. Put
$\nu = a^{1/4}$, and set $z' = z/\nu, u' = - u/(2\nu^2), v' = \nu v$.
Then,
$$\begin{array}{l}
\left[ u',v' \right] = - [u,v]/(2\nu) = z/\nu = z',\\
\left[ u',z' \right] = - [u,z]/(2\nu^3) = \nu v = v',\\
\left[ v',z' \right] = 0.
\end{array}$$
Here we can have the following realization:
$$u' = \left( \begin{array}{rrr}
0 & 1 & 0\\ 1 & 0 & 0\\ 0 & 0 & 0
\end{array} \right),\qquad
v' = \left( \begin{array}{rrr}
0 & 0 & 1\\ 0 & 0 & 0\\ 0 & 0 & 0
\end{array} \right),\qquad
z' = \left( \begin{array}{rrr}
0 & 0 & 0\\ 0 & 0 & 1\\ 0 & 0 & 0
\end{array} \right),$$
and
$$\sigma (u')=u', \sigma(v')=-v', \sigma(z')=-z'.$$
\bigskip

\noindent
{\it Case 2}:\ $c \not= 0,\ a + b = 0$.\par
\noindent
In this case, we have
$[x,y] = z, [z,x] = [z,y] = ax - ay +cz$, and hence
$[u,v] = -2z, [z,u] = 2av + 2cz, [z,v] = 0$.
Set $u' = - u/(2c), v' = cv, z' = z$. Then,
$$\begin{array}{l}
\left[ u',v' \right] = - [u,v]/2 = z = z',\\
\left[ u',z' \right] = - [u,z]/(2c) = av/c + z = av'/(c^2) + z,\\
\left[ v',z' \right] = 0.
\end{array}$$
In this case  we see the following gives a realization for our Lie algebra. 
$$u' = \left( \begin{array}{rrr}
0 & r & 0\\ 1 & 1 & 0\\ 0 & 0 & 0
\end{array} \right),\qquad
v' = \left( \begin{array}{rrr}
0 & 0 & 1\\ 0 & 0 & 0\\ 0 & 0 & 0
\end{array} \right),\qquad
z' = \left( \begin{array}{rrr}
0 & 0 & 0\\ 0 & 0 & 1\\ 0 & 0 & 0
\end{array} \right),$$
and
$$\sigma (u')=u', \sigma(v')=-v', \sigma(z')=-z'$$
with $r = a/(c^2)$.
\bigskip

In order to state a result we put ${\bf P}^1(F) = F \cup \{\ \infty\ \}$. For
each 
$r \in {\bf P}^1(F)$, we define $\Alg(r)$ to be the three dimensional 
Lie algebra spanned by $A(r), X, Y$, where
$$\begin{array}{llllll}
A(r) & = & \left( \begin{array}{rrr}
0 & r & 0\\ 1 & 1 & 0\\ 0 & 0 & 0
\end{array} \right)\quad (r \in F), &
A(\infty) & = & \left( \begin{array}{rrr}
0 & 1 & 0\\ 1 & 0 & 0\\ 0 & 0 & 0
\end{array} \right),\\
&&&&&\\
X & = & \left( \begin{array}{rrr}
0 & 0 & 1\\ 0 & 0 & 0\\ 0 & 0 & 0
\end{array} \right), &
Y & = & \left( \begin{array}{rrr}
0 & 0 & 0\\ 0 & 0 & 1\\ 0 & 0 & 0
\end{array} \right).
\end{array}$$
With this notation we see that the preceding arguments have proven the
following result.
\bigskip

{\bf Theorem 2.5}.\ Let $F$ be an algebraically closed
field of characteristic $0$.
Let $L$ be a three dimensional Lie algebra
with a regular plus-minus pair. Then, $L$ is isomorphic to
one of $sl_2$, ${\mathfrak H}$ or $\Alg(r)$ for $r \in {\bf P}^1(F)$.
\bigskip

We now discuss isomorphisms between these algebras. The dimension of the
derived algebra of $sl_2$ is three, that of  ${\mathfrak H}$ is one while that
of  $\Alg(r)$, for any 
$r \in {\bf P}^1(F)$,  is two so these algebras can never be isomorphic.
Next we
suppose $\Alg(r) \simeq \Alg(\infty)$
with $r \in F$. Then, as is seen by considering the action on the derived
algebra which is spanned by $X,Y$ we see both matrices
$$\left( \begin{array}{rr}
0 & r\\ 1 & 1
\end{array} \right)\qquad \mbox{and}\qquad a
\left( \begin{array}{rr}
0 & 1\\ 1 & 0
\end{array} \right)$$
are similar for some nonzero $a \in F$. Hence,
their characteristic polynomials
$t(t-1) - r$ and $t^2 - a$ must coincide. This however is impossible so
therefore, $\Alg(r) \not\simeq \Alg(\infty)$.
Next we suppose $\Alg(r) \simeq \Alg(s)$ with $r,s \in F$.
Then, both matrices
$$\left( \begin{array}{rr}
0 & r\\ 1 & 1
\end{array} \right)\qquad \mbox{and}\qquad a
\left( \begin{array}{rr}
0 & s\\ 1 & 1
\end{array} \right)$$
are similar for some nonzero $a \in F$. Hence,
their characteristic polynomials
$t(t-1) - r$ and $t(t-a) - as$ must coincide, which implies 
$a = 1$ and $r = s$.
Therefore, $\Alg(r) \not\simeq \Alg(s)$ if $r \not= s$.
\bigskip

Finally we want to remark 
that  the Lie algebra $K(u,1)$ with $u \in F$ 
is isomorphic to
$\Alg(\infty)$ if $u = -1$, and is isomorphic to
$\Alg(- \frac{u}{(u+1)^2})$ if $u \not= -1$. This is because
$$\left( \begin{array}{rr}
-u & 0\ \\ 0 & -1\ 
\end{array} \right)\qquad \mbox{and}\qquad
\left( \begin{array}{rr}
0 & 1\\ 1 & 0
\end{array} \right)$$
are simialr if $u = -1$, and
$$\left( \begin{array}{rr}
-u & 0\ \\ 0 & -1\ 
\end{array} \right)\qquad \mbox{and}\qquad a 
\left( \begin{array}{rr}
0 & r\\ 1 & 1
\end{array} \right)$$
with $a = -(u+1)$ and $r = - \frac{u}{(u+1)^2}$ are similar 
if $u \not= -1$. As a consequence of these computations, we see that
a three dimensional three graded Lie algebra has a regular plus-minus
pair or no plus-minus pair.

\section{\bf Plus-Minus pairs in some general classes of Lie algebras.}

In this section we will show that Borcherds Lie algebras have plus-minus pairs.
Since these generalize the well-known Kac-Moody Lie algebras our results apply
to these as well. We next go on to investigate ${\bf Z}^n$-graded Lie algebras
and see that with certain other assumptions these also have plus-minus pairs.
Our method is to establish slightly more general factorization results in the
universal enveloping algebras of these Lie algebras and then show how this
gives rise to plus-minus pairs. The techniques are quite general and no doubt
apply to other situations as well.
\bigskip

 Let ${\mathfrak g}$ be a rank $l$
Borcherds Lie algebra over $F$ with the  standard
Cartan subalgebra ${\mathfrak h}$ and Chevalley generators
$\{\ e_1,\ldots,e_l,f_1,\ldots,f_l\ \}$, and let ${\mathfrak g}' = 
[ {\mathfrak g},{\mathfrak g} ]$ the derived subalgebra of ${\mathfrak g}$.
Put ${\mathfrak h}' = {\mathfrak h} \cap {\mathfrak g}'$, and take
a compliment ${\mathfrak h}''$ of ${\mathfrak h}'$ in ${\mathfrak h}$
with ${\mathfrak h} = {\mathfrak h}'' \oplus {\mathfrak h}'$. Then,
${\mathfrak g} = {\mathfrak h}'' \oplus {\mathfrak g}'$. Let
${\mathfrak g}_+$ be the subalgebra of ${\mathfrak g}$ generated by
$e_1,\ldots,e_l$, and 
${\mathfrak g}_-$ the subalgebra of ${\mathfrak g}$ generated by
$f_1,\ldots,f_l$ (cf.\ \cite{Bor}, 
\cite{Kac1}, \cite{Kac2}, \cite{Moo}, \cite{MP}, \cite{Wak}).
\bigskip

{\bf Proposition 3.1}.\ Let ${\mathfrak g}$ be a rank $l$ Borcherds Lie
algebra, and let $I \cup J= \{\ 1,2,\ldots,l\ \}$ be a partition of $\{\
1,2,\ldots,l\ \}$ into disjoint subsets. Then,
$$U({\mathfrak g}) = \left( \prod_{i \in I} U(Fe_i) \right)
U({\mathfrak g}_-)\ U({\mathfrak h}'')\ U({\mathfrak g}_+)
\left( \prod_{j \in J} U(Ff_j) \right).$$
\bigskip

{\it Proof }.\ Let ${\mathfrak g}_+^i$ be the standard
homogeneous complementary subalgebra of $Fe_i$ in ${\mathfrak g}_+$, and
${\mathfrak g}_-^i$ the standard
homogeneous complementary subalgebra of $Ff_i$ in ${\mathfrak g}_-$.
For each $k = 1,\ldots,l$ we put $h_k = [e_k,f_k]$ and
${\mathfrak h}_k = Fh_{k+1} \oplus \cdots \oplus Fh_l$, and we set
$I_k = I \cap \{\ 1,\ldots,k\ \}$ and $J_k = J \cap \{\ 1,\ldots,k\ \}$. We
make
free use of the
PBW Theorem as well as the fact 
that $Fe_i \oplus Ff_i \oplus Fh_i$ is either $sl_2$ or $\mathfrak H$ so has a
regular plus-minus pair.

If $1 \in I$, then
$$\begin{array}{lll}
U({\mathfrak g}) & = & U({\mathfrak g}_-)\ U({\mathfrak h})\ 
U({\mathfrak g}_+)\\
& = & U({\mathfrak g}_-^1)\ U(Ff_1)\ U({\mathfrak h}'')\ 
U({\mathfrak h}_1)\ U(Fh_1)\ U(Fe_1)\ U({\mathfrak g}_+^1)\\
& = & U({\mathfrak g}_-^1)\ U({\mathfrak h}'')\ 
U({\mathfrak h}_1)\ U(Ff_1)\ U(Fh_1)\ U(Fe_1)\ 
U({\mathfrak g}_+^1)\\
& = & U({\mathfrak g}_-^1)\ U({\mathfrak h}'')\ U({\mathfrak h}_1)\
U(Fe_1)\ U(Ff_1)\ U(Fe_1)\ U({\mathfrak g}_+^1)\\
& = & U(Fe_1)\ U({\mathfrak g}_-^1)\ U({\mathfrak h}'')\ U({\mathfrak h}_1)\ 
U(Ff_1)\ U(Fe_1)\ U({\mathfrak g}_+^1)\\
& = & U(Fe_1)\ U({\mathfrak g}_-)\ U({\mathfrak h}'')\ U({\mathfrak h}_1)\ 
U({\mathfrak g}_+).
\end{array}$$
In the other case when $1 \in J$ by using the same type of argument we have
$$U({\mathfrak g}) =  
U({\mathfrak g}_-)\ U({\mathfrak h}'')\ U({\mathfrak h}_1)\
U({\mathfrak g}_+)\ U(Ff_1).$$
If we began with
$$U({\mathfrak g}) = \left( \prod_{i \in I_k} U(Fe_i) \right)
U({\mathfrak g}_-)\ U({\mathfrak h}'')\ U({\mathfrak h}_k)\ 
U({\mathfrak g}_+) \left( \prod_{j \in J_k} U(Ff_j) \right),$$
then, again using the same method, we can obtain
$$U({\mathfrak g}) = \left( \prod_{i \in I_{k+1}} U(Fe_i) \right)
U({\mathfrak g}_-)\ U({\mathfrak h}'')\ U({\mathfrak h}_{k+1})\
U({\mathfrak g}_+) \left( \prod_{j \in J_{k+1}} U(Ff_j) \right).$$
Thus after several applications of this process we reach the stated result. 
Q.E.D.
\bigskip

We next see that this gives the desired plus-minus pair.
\hfill
\eject

{\bf Corollary 3.2}.\ Let $\mathfrak g$ be a Borcherds Lie algebra of finite
rank. Then,
$$U({\mathfrak g}) = U({\mathfrak g}_\pm)\ 
U({\mathfrak g}_\mp)\ U({\mathfrak h}'')\ U({\mathfrak g}_\pm).$$
Hence,  perfect Borcherds Lie algebras have plus-minus pairs.
In particular, perfect Kac-Moody Lie algebras or, more generally, perfect
Brocherds Lie have regular plus-minus
pairs.
\bigskip

{\it Proof }.\ We just take 
one of $I$ and $J$ to be empty. This leads to the result. 
Then, for example, let $P = {\mathfrak h}'' \oplus {\mathfrak g}_+$ and
$M = {\mathfrak g}_-$. This gives a plus-minus pair. Q.E.D.
\bigskip

We next generalize the previous discussion by considering ${\bf
Z}^n$-graded Lie
algebras. Thus, let $Q = \oplus_{i = 1}^{n} {\bf Z} \alpha_i$ 
be a free abelian
group of rank $n$ generated by $\alpha_1,\ldots,\alpha_n$, and let ${\mathfrak
g}
= \oplus_{\alpha \in Q} {\mathfrak g}_\alpha$ be a Lie algebra graded by $Q$.
Put $\Delta = \{\ \alpha \in Q\ 
\mid\ {\mathfrak g}_\alpha \not= 0\ \}$.
We also assume  that ${\bf Z}\alpha_1 \cap 
\Delta = \{\ 0,\ \pm \alpha_1\ \}$, and that $L = 
{\mathfrak g}_{\alpha_1} \oplus [{\mathfrak g}_{\alpha_1},
{\mathfrak g}_{- \alpha_1}] \oplus {\mathfrak g}_{- \alpha_1}$
is a subalgebra with
a plus-minus pair $({\mathfrak g}_{\alpha_1},{\mathfrak g}_{- \alpha_1})$ 
in $L$.
(Thus, if $L$ is
three dimensional Proposition 2.4 implies $L$ is isomorphic to 
either $sl_2$ or $\mathfrak H.$) We also suppose
that there exists a complementary subalgebra ${\mathfrak g}_0'$ of 
$[{\mathfrak g}_{\alpha_1},{\mathfrak g}_{- \alpha_1}]$ in ${\mathfrak g}_0$
with ${\mathfrak g}_0 =
[{\mathfrak g}_{\alpha_1},{\mathfrak g}_{- \alpha_1}] \oplus 
{\mathfrak g}_0'$. An element $\alpha = \sum_{i=1}^{n} c_i\alpha_i \in Q$
is called positive (resp. negative), that is $\alpha > 0$, if
there is an index $i$ satisfying $c_i > 0$ (resp. $c_i < 0$) and
$c_{i+1} = c_{i+2} = \ldots = c_{n} = 0$. Put $\Delta_+ = 
\{\ \alpha \in \Delta\ \mid\ \alpha > 0\ \}$ and $\Delta_- = 
\{\ \alpha \in \Delta\ \mid\ \alpha < 0\ \}$. Let ${\mathfrak g}_\pm
= \oplus_{\alpha \in \Delta_\pm} {\mathfrak g}_\alpha$, and
${\mathfrak g}_\pm'
= \oplus_{\alpha \in \Delta_\pm \setminus \{ \alpha_1 \}} 
{\mathfrak g}_\alpha$. Then, ${\mathfrak g}_\pm = 
{\mathfrak g}_{\pm \alpha_1} \oplus {\mathfrak g}_\pm'$, and we see that
${\mathfrak g}_{\pm \alpha_1} \oplus {\mathfrak g}_\mp'$ are subalgebras.
In this situation we have the following result.
\bigskip

{\bf Proposition 3.3}.\ Let ${\mathfrak g} = 
\oplus_{\alpha \in Q} {\mathfrak g}_\alpha$
be a graded Lie algebra with the extra conditions  as above. Then,
$$\begin{array}{lll}
U({\mathfrak g}) & = & U({\mathfrak g}_{\alpha_1})\ 
U({\mathfrak g}_-)\ U({\mathfrak g}_0')\ U({\mathfrak g}_+)
\end{array}$$ Moreover, letting $P = {\mathfrak g}_+ \oplus
{\mathfrak g}_0'$ and $M = {\mathfrak g}_-$ gives a plus-minus pair for
$\mathfrak g$.

\bigskip

{\it Proof. }  Using our assumptions we see that
$$\begin{array}{lll}
U({\mathfrak g}) & = & U({\mathfrak g}_-)\ U({\mathfrak g}_0)\ 
U({\mathfrak g}_+)\\
& = & U({\mathfrak g}_-')\ U({\mathfrak g}_{- \alpha_1})\ 
U({\mathfrak g}_0')\ 
U([ {\mathfrak g}_{\alpha_1} , {\mathfrak g}_{- \alpha_1} ])\ 
U({\mathfrak g}_{\alpha_1})\ U({\mathfrak g}_+')\\
& = & U({\mathfrak g}_-')\ U({\mathfrak g}_0')\ 
U({\mathfrak g}_{- \alpha_1})\ 
U([ {\mathfrak g}_{\alpha_1} , {\mathfrak g}_{- \alpha_1} ])\ 
U({\mathfrak g}_{\alpha_1})\ U({\mathfrak g}_+')\\
& = & U({\mathfrak g}_-')\ U({\mathfrak g}_0')\ 
U({\mathfrak g}_{\alpha_1})\ U({\mathfrak g}_{- \alpha_1})\ 
U({\mathfrak g}_{\alpha_1})\ U({\mathfrak g}_+')\\
& = & U({\mathfrak g}_{\alpha_1})\ U({\mathfrak g}_-')\ 
U({\mathfrak g}_0')\ U({\mathfrak g}_{- \alpha_1})\ 
U({\mathfrak g}_{\alpha_1})\ U({\mathfrak g}_+')\\
& = & U({\mathfrak g}_{\alpha_1})\ U({\mathfrak g}_-')\ 
U({\mathfrak g}_{- \alpha_1})\ U({\mathfrak g}_0')\ 
U({\mathfrak g}_{\alpha_1})\ U({\mathfrak g}_+')\\
& = & U({\mathfrak g}_{\alpha_1})\ U({\mathfrak g}_-)\ 
U({\mathfrak g}_0')\ U({\mathfrak g}_+).
\end{array}$$
Q.E.D.
\bigskip

{\bf Remark}\quad It can be seen that many EALA's and some of the root-graded
Lie algebras (cf.\ 
\cite{AABGP}, \cite{ABG}, \cite{BM}) satisfy the hypothesis of Proposition 3.3.
\bigskip



\begin{thebibliography}{99}




\bibitem{AABGP} B.\ N.\ Allison, S.\ Azam, S.\ Berman, Y.\ Gao and A.\
Pianzola , Extended affine Lie algebras and their root systems, Memoirs
AMS 126 (No.\ 603), Providence, 1997.




\bibitem{ABG} B.\ N.\ Allison, G.\ Benkart and Y.\ Gao , 
Central extensions of Lie algebras graded by finite root systems,
Math.\ Ann.\ 316 (2000), 499 -- 527.




\bibitem{BM} S.\ Berman and R.\ V.\ Moody , Lie algebras graded by finite 
root systems and the intersection matrix algebras of Slodowy, Invent.\ Math.\ 
108 (1992), 323 -- 347.




\bibitem{Bor} R.\ E.\ Borcherds , Generalized Kac-Moody algebras,\
J.\ Algebra 115 (1988), 501 -- 512.




\bibitem{Hum} J.\ E.\ Humphreys , ``Introduction to Lie algebras and
representation theory'', GTM 9, Springer-Verlag, New York, 1972.




\bibitem{Jac} N.\ Jacobson , ``Lie algebras'', Interscience, New York, 1962.




\bibitem{Kac1}\ V.\ G.\ Kac , Simple irreducible graded Lie algebras of finite
growth, Math.\ USSR-Izv.\ 2 (1968), 211 -- 230.




\bibitem{Kac2} V.\ G.\ Kac , ``Infinite dimensional Lie algebras'',
3rd edition, Cambridge University Press, Cambridge, 1990.




\bibitem{Moo} R.\ V.\ Moody , 
A new class of Lie algebras, J.\ Algebra 10 (1968), 211 -- 230.




\bibitem{MP} R.\ V.\ Moody and A.\ Pianzola , ``Lie algebras with triangular
decompositions'', J.\ Wiley \& Sons, New York, 1995.




\bibitem{Tits} J.\ Tits , Th\' eorie des groupes, ``R\' edum\' e des cours et 
travaux (1980 -- 1981)'', pp.\ 75 -- 87, Coll\' ege de France, Paris, 1981.




\bibitem{Wak} M.\ Wakimoto , ``Infinite-dimensional Lie algebras'',
Translations of Mathematical Monographs, Iwanami Series of Modern Mathematics,
2001.




\end{thebibliography}
\end{document}